\documentclass[10pt]{amsart}
\usepackage{graphicx}
\usepackage{tikz}
\newtheorem{thm}{Theorem}[section]
\newtheorem{cor}[thm]{Corollary}
\newtheorem{lem}[thm]{Lemma}

\theoremstyle{definition}

\theoremstyle{remark}

\numberwithin{equation}{section}

\title{The Critical Independence Number and an Independence Decomposition}
\author{C. E. Larson\\Department of Mathematics and Applied Mathematics\\Virginia Commonwealth University\\Richmond, VA 23284\\clarson@vcu.edu\\804-828-5576}

\begin{document}

\begin{abstract}
An independent set $I_c$ is a \textit{critical independent set} if $|I_c| - |N(I_c)| \geq |J| - |N(J)|$, for any independent set $J$. 
The \textit{critical independence number} of a graph is the cardinality of a maximum critical independent set. This number is a lower bound for the independence number and can be computed in polynomial-time.  Any graph can be decomposed into two subgraphs where the independence number of one subgraph equals its critical independence number, where the critical independence number of the other subgraph is zero, and where the sum of the independence numbers of the subgraphs is the independence number of the graph.  A proof of a conjecture of Graffiti.pc yields a new characterization of K\"{o}nig-Egervary graphs: these are exactly the graphs whose independence and critical independence numbers are equal.
\end{abstract}

\maketitle

\section{Introduction}

An independent set of vertices in a graph is a set of vertices no two of which are adjacent. A maximum independent set is an independent set of largest cardinality. Finding a maximum independent set (MIS) in a graph is a well-known widely-studied NP-hard problem \cite{GareJohn}. It will be shown that the problem of finding a MIS in a graph $G$ can be decomposed into finding a MIS for two subgraphs, $G[X]$ and $G[X^c]$, where $X$ is a maximum critical independent set together with its neighbors, and $X^c=V(G)\setminus X$. The union of these  independent sets is a MIS in $G$. There is an efficient algorithm for finding both the set $X$ and a MIS in $G[X]$.

The following notation is used throughout: the vertex set of a graph $G$ is $V(G)$, the
order of $G$ is $n=n(G)=|V(G)|$, the set of neighbors of a vertex $v$ is $N_G(v)$ (or simply $N(v)$ if there is no possibility of ambiguity), the set of neighbors of a set $S \subseteq V(G)$ in $G$ is $N_G(S)=\cup_{u\in S}N(u)$ (or simply $N(S)$ if there is no possibility of ambiguity), the graph induced on $S$ is $G[S]$, and the independence number, the cardinality of a MIS, is $\alpha(G)$. All graphs are assumed to be finite
and simple.

An independent set of vertices $I_c$ is a \emph{critical independent set} if $|I_c| - |N(I_c)|\geq |J| - |N(J)|$, for any independent set $J$.   A graph may contain critical independent sets of different cardinalities. A graph consisting of a single edge ($K_2$, the complete graph on two vertices) has critical independent sets of cardinalities $0$ and $1$. For some graphs the only critical independent set is the empty set; $K_3$ is an example.
A \emph{maximum critical independent set} is a critical independent set of maximum cardinality. It is easy to verify that, for any graph with at least three vertices, a maximum critical independent set must contain all pendant vertices; so a maximum critical independent set is a generalization of the set of pendants. The \emph{critical independence number} of a graph $G$, denoted $\alpha'=\alpha'(G)$, is the cardinality of a maximum critical independent set. If $I_c$ is a maximum critical independent set, and so $\alpha'(G)=|I_c|$, then clearly $\alpha'\leq \alpha$.  
Much of the interest in critical independent sets is due to the following theorem.
\begin{thm} \label{t:butenko} (Butenko \& Trukhanov, \cite{ButeTruk}) 
If $I_c$ is a critical independent set in a graph $G$ then there is a maximum independent set $I$ in $G$ such that $I_c\subseteq I$.
\end{thm} 
Butenko and Trukhanov also proposed the problem of finding a polynomial-time algorithm for finding a maximum critical independent set in a graph \cite{ButeTruk}. Their problem was solved by this author \cite{Lars};
thus the critical independence number of a graph can be computed in polynomial-time.

\begin{figure}
\begin{center}
\includegraphics{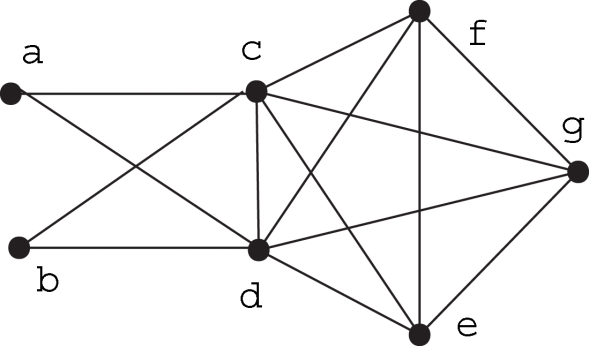}
\caption{The vertices $I_c=\{a,b\}$ form a (maximum cardinality) critical independent set; according to Theorem \ref{t:butenko}, this set of vertices can be extended to a maximum independent set of the graph.  The sets $X=I_c\cup N(I_c)=\{a,b,c,d\}$ and $X^c=V\setminus X=\{e,f,g\}$ induce a decomposition of the graph into a totally independence reducible subgraph $G[X]$ and an independence irreducible subgraph $G[X^c]$, according to Theorem \ref{t:decomposition}.}
\end{center}
\label{fig:critical-independent}
\end{figure}

A graph   is \emph{independence irreducible} if $\alpha'=0$. This means that the empty set is the only critical independent set. It is easy to see that a graph is independence irreducible if, and only if, the number of neighbors of any non-empty independent set of vertices is greater than the number of vertices in the set. Complete graphs with at least three vertices, odd cycles, and fullerene graphs \cite{Lars} are examples.  A graph is \emph{independence reducible} if $\alpha'>0$. This means that the graph is guaranteed to have a non-empty critical independent set. A graph is \emph{totally independence reducible} if $\alpha'=\alpha$; $K_2$ and even cycles are examples. 
\begin{figure}
\begin{center}
\includegraphics[width=4in]{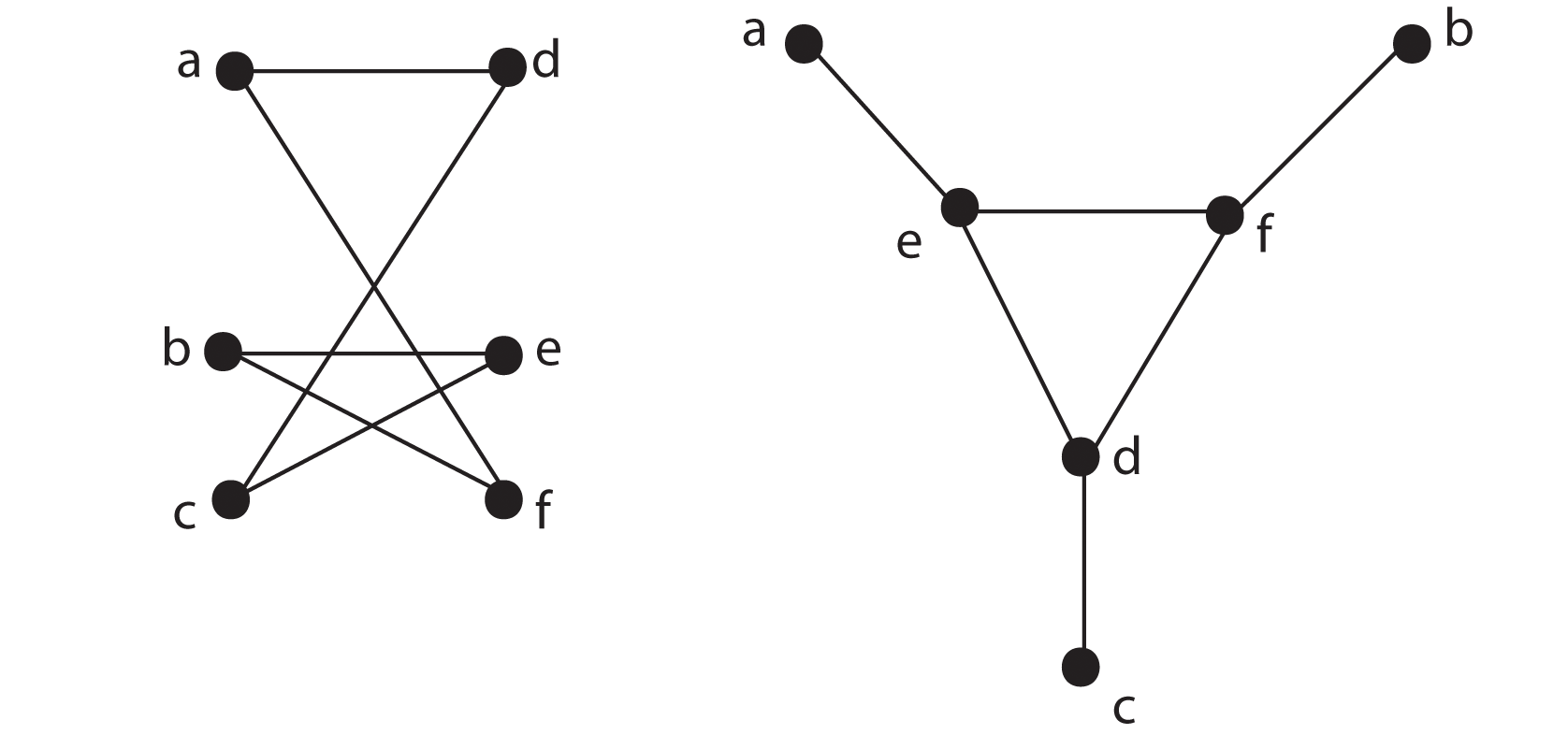}
\caption{These graphs are \textit{totally independence reducible}. In the  graph on the left, the set $I=\{a,b,c\}$ is a maximum critical independent set, and a maximum independent set; all bipartite graphs are totally independence reducible. The graph on the right is not bipartite; the set $I_c=\{a,b,c\}$ is a maximum critical independent set, and a maximum independent set. For both of these graphs $\alpha=\alpha'=3$.}
\end{center}
\label{f:2}
\end{figure}

\begin{figure}
\begin{center}
\includegraphics[width=5in]{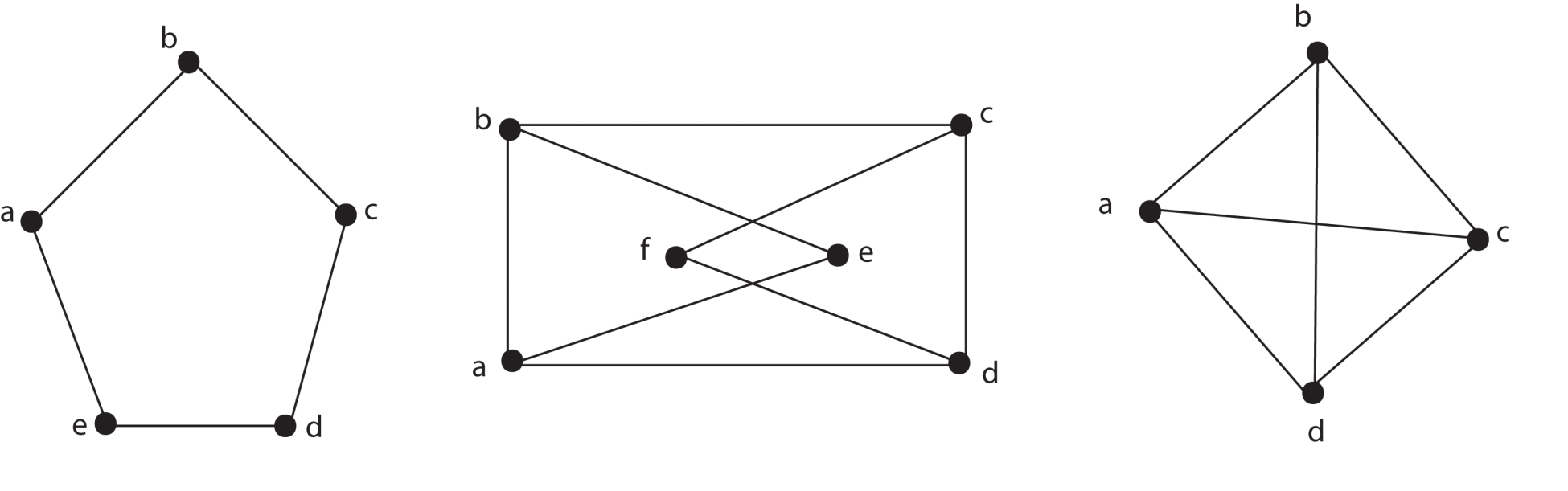}\label{fig:independence_irreducible}
\caption{These graphs are independence irreducible: for \textit{any} non-empty independent set of vertices $I$ in these graphs, $|N(I)|>|I|$.}
\end{center}
\label{fig:irreducible}
\end{figure}

The definition of a critical independent set is due to Zhang's 1990 paper \cite{Zhan}; he showed that these could be found in polynomial time. In 1994 Ageev \cite{Agee} provided a simpler algorithm, reducing the problem to that of finding a maximum independent set in a bipartite graph.  Then, after more than ten years elapsed,  Butenko and Trukhanov \cite{ButeTruk} proved their Theorem \ref{t:butenko}, thereby directly connecting the problem of finding a critical independent set to that of finding a maximum independent set.

\section{An Independence Decomposition}\label{s:decomposition}

Finding a maximum independent set in a graph $G$ and computing its independence number are NP-hard problems. When attacking these problems it would be useful to be able to decompose the problem into finding maximum independent sets for the graphs induced by the sets in some partition of the vertex set $V(G)$. It will be shown that a non-trivial partition exists and, furthermore, that an efficient algorithm exists for finding a MIS of at least one of the corresponding subgraphs.

Butenko and Trukhanov \cite{ButeTruk} noted that, if the independence number of a graph is at least half the number of vertices then the graph will have a non-empty critical independent set; the idea is that either a maximum independent set $I$ will be a critical independent set, since $|I| - |N(I)| \geq 0$, or there must be a non-empty independent set $J$ such that $|J| - |N(J)| > |I| - |N(I)|$. In either case $\alpha'>0$, and the graph is independence reducible. Furthermore, this means that independence irreducible graphs have ``small'' (less than $\frac{n}{2}$) independence numbers.

The following characterization of graphs whose independence numbers equal their critical independence numbers will be needed in the proof of the main result.

\begin{lem}\label{t:char}
For any graph $G$ with maximum critical independent set $I_c$, $\alpha=\alpha'$ if, and only if, $I_c \cup N(I_c)=V(G)$.
\end{lem}

\begin{proof}
Let $G$ be a graph. Suppose first that $\alpha(G)=\alpha'(G)$. Let $I_c$ be a maximum critical independent set of $G$. So $\alpha'=|I_c|$. Suppose $I_c\cup N(I_c)$ is a proper subset of $V(G)$. Let $v\in V\setminus (I_c\cup N(I_c))$. Then $I_c\cup \{v\}$ is an independent set and $\alpha(G)\geq |I_c|+1=\alpha'(G)+1>\alpha'(G)$, contradicting the fact that $\alpha(G)=\alpha'(G)$.

Suppose now that $I_c$ is a maximum critical independent set and $I_c\cup N(I_c)=V(G)$. Theorem \ref{t:butenko} implies that there is a maximum independent set $I$ of $G$, such that $I_c \subseteq I$. If there is a vertex $v\in I\setminus I_c$ then, by assumption, $v\in N(I_c)$. But then $v$ is adjacent to some vertex in $I_c$ and $I$ is not independent. So $I=I_c$ and $\alpha=\alpha'$.
\end{proof}

Since a maximum critical independent set of a graph can be found in polynomial-time, Lemma \ref{t:char} implies that whether a graph has the property that $\alpha=\alpha'$ (that is, whether the graph is totally independence irreducible) can be determined in polynomial-time.

A \textit{matching} in a graph is a set of non-incident edges. The \textit{matching number} $\mu$ is the cardinality of a largest matching. A matching $M$  is a \textit{matching of a set $X$ into a disjoint set $Y$} if every vertex in $X$ is incident to some edge in $M$ and each of these edges is incident to a vertex in $Y$; it is not required that every vertex in $Y$ be incident to an edge in $M$. 
\begin{lem}\label{t:matching} (The Matching Lemma, Larson, \cite{Lars}) If $I_c$ is a critical independent set of $G$, then there is a matching of the vertices $N(I_c)$ into (a subset of) the vertices of $I_c$.
\end{lem}
The proof of the Matching Lemma is essentially an application of Hall's Theorem. 

\begin{lem}\label{l:crit}
If $G$ is a graph with critical independent sets $I_c$ and $J_c$, where $J=J_c\setminus (I_c\cup N(I_c))$, and $I=I_c\cup J$ then,
\begin{enumerate}
\item $|I_c\cap N(J_c)|=|J_c\cap N(I_c)|$,
\item $|J|\geq |N(J_c)\setminus (I_c\cup N(I_c)|$, and
\item $I$ is a critical independent set.
\end{enumerate}
\end{lem}

\begin{figure}
\begin{tikzpicture}[scale=1]
\shade (0,0) rectangle (3,7);
\shade (6,0) rectangle (10,2);
\foreach \x in {0,3,6,10} \foreach \y in {0,2,4} \draw (\x,\y) circle (1.5pt) [fill=black];
\foreach \x in {0,3,6} \foreach \y in {7} \draw (\x,\y) circle (1.5pt) [fill=black];
\foreach \x in {0,3,6,10} \draw  [line width=1pt] (\x,0)--(\x,2);
\foreach \x in {0,3,6,10} \draw  [line width=1pt] (\x,4)--(\x,2);
\foreach \x in {0,3,6} \draw  [line width=1pt] (\x,7)--(\x,4);
\foreach \y in {0,2,4,7} \draw  [line width=1pt] (0,\y)--(3,\y);
\foreach \y in {0,2,4,7} \draw  [line width=1pt] (3,\y)--(6,\y);
\foreach \y in {0,2,4} \draw  [line width=1pt] (6,\y)--(10,\y);
\draw (1.5,7.4) node {$I_c$}; 
\draw (4.5,7.4) node {$N(I_c)$}; 
\draw (10.7,3) node {$N(J_c)$}; 
\draw (10.5,1) node {$J_c$}; 

\draw (1.5,1) node {$I_c\cap J_c$}; 
\draw (1.5,3) node {$I_c\cap N(J_c)$}; 
\draw (4.5,3) node {$N(I_c)\cap N(J_c)$}; 
\draw (4.5,1) node {$N(I_c)\cap J_c$}; 
\draw (8,1) node {$J_c\setminus (I_c\cup N(I_c))$}; 

\end{tikzpicture}

\caption{A useful figure for following the proofs of Lemma \ref{l:crit} and Theorem \ref{t:decomposition}. The figure  is a schematic of the relationship between critical independent sets $I_c$ and $J_c$ and their neighbors.  The  set $I=I_c\cup J=I_c\cup J_c \setminus(I_c \cup N(I_c))$ is shaded.}
\label{fig:criticalsets}
\end{figure}
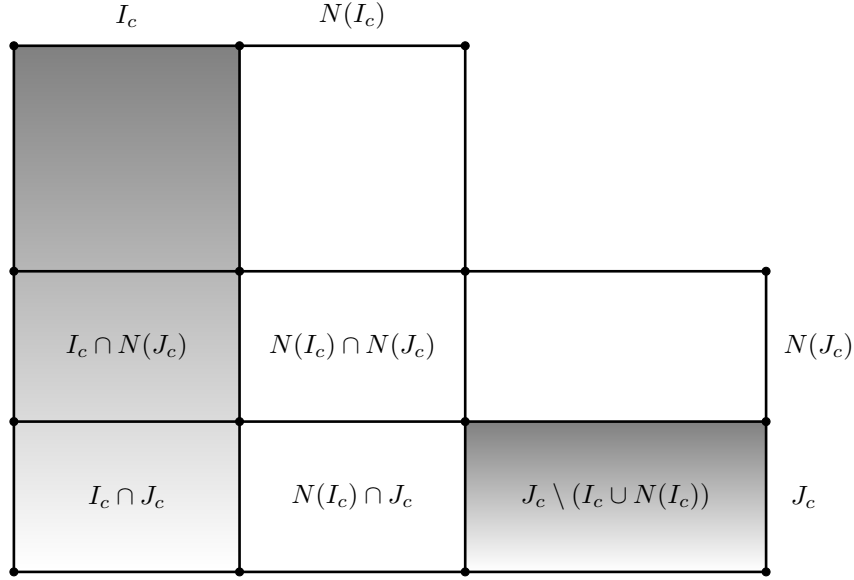

\begin{proof}
The Matching Lemma \ref{t:matching} guarantees that there is a matching from the vertices in $N(J_c)$ to (a subset of) the vertices in $J_c$ and
from the vertices in $N(I_c)$ to (a subset of) the vertices in $I_c$. For the remainder of the proof the reader may usefully refer to Figure \ref{fig:criticalsets}. Since the vertices in $I_c\cap N(J_c) \subseteq N(J_c)$ must be matched to vertices in $N(I_c)\cap J_c$, and the vertices in
$N(I_c)\cap J_c \subseteq N(I_c)$ must be matched to vertices in $I_c \cap N(J_c)$, it follows that $|I_c\cap N(J_c)|=|J_c\cap N(I_c)|$, proving (1).

Applying the Matching Lemma again, we have that $N(J_c)$ is matched into $J_c$, that is, every vertex in $N(J_c)$ can be paired with a distinct adjacent vertex in $J_c$. Notice that a vertex $v$ in $N(J_c)\setminus (I_c\cup N(I_c))$ cannot be matched to a vertex in $J_c \cap N(I_c)$ under any matching, as the proof of (1) guarantees that these are only matched to vertices in $I_c \cap N(J_c)$. Furthermore, a vertex $v$ in $N(J_c)\setminus (I_c\cup N(I_c))$ cannot be matched to a vertex $w$ in $I_c \cap J_c$. If it were, then since $w \in I_c$ and $v$ is adjacent to $w$, it follows that $v\in N(I_c)$, contradicting the fact that $v \notin N(I_c)$. Thus vertices in $N(J_c)\setminus (I_c\cup N(I_c))$ can only be matched to vertices in $J_c \setminus(I_c \cup N(I_c))$. Since every vertex in $N(J_c)\setminus (I_c\cup N(I_c))$ \textit{is} matched to a vertex in $J_c \setminus(I_c \cup N(I_c))$, it follows that $|J|=|J_c \setminus(I_c \cup N(I_c))|\geq |N(J_c)\setminus (I_c\cup N(I_c))|$, proving (2).

$I=I_c\cup J$. Since $I_c$ and $J$ are independent, and $J=J_c\setminus (I_c\cup N(I_c))$, $I$ is independent.
Since $I_c$ and $J$
are disjoint, $|I| = |I_c| + |J|$. $N(I)\subseteq N(I_c)\cup [N(J_c)\setminus (I_c\cup N(I_c))]$ and
$|N(I)| \leq |N(I_c)| + |N(J_c) \setminus (I_c\cup N(I_c))|$.
So,
$
|I| - |N(I)| \geq (|I_c| + |J|) - (|N(I_c)| + |N(J_c)\setminus (I_c\cup N(I_c))| = (|I_c| - |N(I_c)|) + (|J| - |N(J_c)\setminus (I_c\cup N(I_c))|).
$
Since (2) implies that the last term is non-negative, it follows that $|I|-|N(I)|\geq |I_c|-|N(I_c)|$ and, thus, that $I$ is a critical independent set, proving (3).
\end{proof}

\begin{thm}\label{t:decomposition}
For any graph $G$, there is a unique set $X\subseteq V(G)$ such that
\begin{enumerate}
\item $\alpha(G)=\alpha(G[X])+\alpha(G[X^c])$,
\item $G[X]$ is totally independence reducible,
\item $G[X^c]$ is independence irreducible, and
\item for every maximum critical independent set $J_c$ of $G$, $X=J_c\cup N(J_c)$.
\end{enumerate}
\end{thm}

\begin{proof}
Let $I_c$ be a maximum critical independent set of $G$. Let $X=I_c\cup N(I_c)$ and $X^c=V(G)\setminus X$.  $I_c$ is an independent set in $G[X]$.

 Suppose $I_c$ is not a maximum independent set in $G[X]$. Let $Y$ be an independent set of $G[X]$ such that $|Y|>|I_c|$. Let $Y_I=Y\cap I_c$ and $Y_N=Y\cap N(I_c)$. So $Y=Y_I\cup Y_N$,  $|Y_I| + |Y_N|=|Y|$, and $|Y_I| > |I_c| - |Y_N|$. Note that $N(Y_I)\subseteq N(I_c)\setminus Y_N$. Then, $|Y_I| - |N(Y_I)| \geq |Y_I| - |N(I_c)\setminus Y_N| > (|I_c| - |Y_N|) - |N(I_c)\setminus Y_N| = |I_c| - (|Y_N| + |N(I_c)\setminus Y_N|)=|I_c| - |N(I_c)|$. Since $Y_I$ is an independent set, $I_c$ cannot be a critical independent set of $G$, which contradicts the assumption that it is. Thus, $I_c$ \textit{is} a maximum independent set and $\alpha(G[X])=|I_c|$.

It follows from Theorem \ref{t:butenko} that $I_c$ is contained in a maximum independent set $I$ of $G$. So $\alpha(G)=|I|$. $I\setminus I_c$ is an independent set in $X^c$. So $\alpha(G[X^c])\geq |I\setminus I_c|$. Suppose there is an independent set $I'\subseteq X^c$ such that $|I'|>|I\setminus I_c|$. By construction, no vertex in $I_c$ is adjacent in $G$ to a vertex in $X^c$. Thus, no vertex in $I_c$ is adjacent to a vertex in $I'$. Thus,  $I_c\cup I'$ is an independent set in $G$, and $\alpha(G)\geq |I_c\cup I'|=|I_c|+|I'|>|I_c|+|I\setminus I_c|=|I|=\alpha(G)$, a contradiction. Thus,  $I\setminus I_c$ is a maximum independent set in $G[X^c]$, and $\alpha(G[X])+\alpha(G[X^c])=|I_c|+|I\setminus I_c|=|I|=\alpha(G)$, proving (1).

Now suppose $I_c$ is not a critical independent set in $G[X]$. Let $Y$ be a minimum critical independent set of $G[X]$.  So  $|Y|-|N_{G[X]}(Y)|>|I_c|-|N(I_c)|$. Let $Y_I=Y\cap I_c$ and $Y_N=Y\cap N(I_c)$. (Note that $N(I_c)$ is unambiguous as $N_G(I_c)=N_{G[X]}(I_c)$.)
Let $Y_N'\subseteq I_c$ be the set of neighbors of $Y_N$ in $I_c$. It follows from the Matching Lemma \ref{t:matching} that there is a matching of the vertices in $N_{G[X]}(Y)$ to (a subset of) the vertices in $Y$. Since $I_c$ is an independent set, and
$Y_N' \subseteq I_c$, the vertices in $Y_N'$ must be matched to vertices in $Y_N$. Thus, $|Y_N|\geq |Y_N'|$.

Suppose $|Y_N|=|Y_N'|$. Then $|Y_I|-|N(Y_I)|=(|Y_I|+|Y_N|)-(|N(Y_I)|+|Y_N|)=|Y|-(|N(Y_I)|+|Y_N'|)\geq |Y|-|N_{G[X]}(Y)|$, implying that $Y_I$ is a critical independent set of $G[X]$. Since $Y_I\subseteq Y$, and $Y$ is a minimum critical independent set, it follows that $Y_I=Y$, and $Y_N=\emptyset$. Since $Y_I\subseteq I_c$, $N_G(Y_I)=N_{G[X]}(Y_I)$, and $|Y_I|-|N_G(Y_I)|\geq |Y|-N_{G[X]}(Y)|>|I_c|-|N_G(I_c)|$, contradicting the fact that $I_c$ is a critical independent set in $G$.

So $|Y_N|>|Y_N'|$.  But then, for $I=I_c\setminus Y_N'$, $|I|-|N_G(I)|=|I_c\setminus Y_N'|-|N(I_c)\setminus Y_N|=|I_c|-|Y_N'|-(|N(I_c)|-|Y_N|)=(|I_c|-|N(I_c)|)+(|Y_N|-|Y_N'|)$. Since the last term is positive, it follows that $|I|-|N_G(I)|>|I_c|-|N(I_c)|$, again contradicting the fact that $I_c$ is a critical independent set in $G$. Thus, $I_c$ is a critical independent set in $G[X]$. Since $I_c\cup N(I_c)=X$, $I_c$ is a maximum critical independent set in $G[X]$, and $\alpha'(G[X])=|I_c|$.
So $\alpha(G[X])=\alpha'(G[X])=|I_c|$ and $G[X]$ is totally independence reducible, proving (2).

Suppose that $G[X^c]$ contains a non-empty critical independent set $Z$. So $|Z|\geq |N_{G[X^c]}(Z)|$.  No vertex in $I_c$ is  adjacent to any vertex in $Z$ as $N(I_c)\subseteq X$ and $Z\subseteq X^c$. So $I_c\cup Z$ is an independent set in $G$. Furthermore, $|N(I_c\cup Z)|=|N(I_c)|+|N_{G[X^c]}(Z)|$. So, $|I_c\cup Z|-|N_G(I_c\cup Z)|= (|I_c|+|Z|)-(|N(I_c)|+|N_{G[X^c]}(Z)|)=(|I_c|-|N(I_c)|)+(|Z|-|N_{G[X^c]}(Z)|)\geq|I_c|-|N(I_c)|$, contradicting the fact that $I_c$ is a maximum critical independent set of $G$. Thus, $G[X^c]$ does not contain a non-empty critical independent set, $\alpha'(G[X^c])=0$, and $G[X^c]$ is irreducible, proving (3).

Now suppose that $J_c$ is a maximum critical independent set of $G$. Thus, since $J_c$ and $I_c$ are both maximum critical independent sets, $|J_c|=|I_c|$. Since they are both critical, $|J_c|-|N(J_c)|=|I_c|-|N(I_c)|$. It then follows that $|N(J_c)|=|N(I_c)|$.  Let $J=J_c\setminus (I_c\cup N(I_c))$. So $I_c\cup J$ is an independent set. Lemma  \ref{l:crit} implies that $I_c\cup J$ is a critical independent set of $G$. But, since $I_c\subseteq I_c\cup J$ and $I_c$ is a maximum critical independent set of $G$, $J=\emptyset$.  A parallel argument yields that $I=\emptyset$.

The Matching Lemma \ref{t:matching} implies that there is a matching from the vertices in $N(J_c)$ into the vertices in $J_c$. Lemma \ref{l:crit} implies that $|I_c\cap N(J_c)|=|J_c\cap N(I_c)|$. So if $v\in N(J_c)\setminus (N(J_c)\cap I_c)$, it must be matched to a vertex in $J_c\setminus (J_c\cap N(I_c)=I_c\cap J_c$ and, thus, $v\in N(I_c\cap J_c)\subseteq N(I_c)\cap N(J_c)$. So every vertex in $N(J_c)$ is either in $N(J_c)\cap I_c$ or in $N(J_c)\cap N(I_c)$, which implies that $N(J_c)\subseteq I_c\cup N(I_c)$.

So both $J_c$ and $N(J_c)$ are subsets of $I_c\cup N(I_c)$. Since $|J_c|+|N(J_c)|=|I_c|+|N(I_c)|$, it follows that $J_c\cup N(J_c)=I_c \cup N(I_c)=X$, proving (4).

The uniqueness of a set $X\subseteq V(G)$ satisfying the four conditions of the theorem follows immediately from (4).

\end{proof}

\section{An application: K\"{o}nig-Egervary Graphs}
\label{s:konig}

The independence number, the critical independence number, order, and the matching number of a graph are  $\alpha$,  $\alpha'$, $n$ and $\mu$, respectively. A \textit{vertex cover} in a graph is a set of vertices such that each edge in the graph is incident to at least one of the vertices in the cover. The \textit{vertex covering number} $\tau$ is the cardinality of a smallest vertex cover. One of the Gallai Identities is that, for any  graph, $\alpha+\tau=n$ \cite[p.~2]{LovaPlum}. For bipartite graphs, $\alpha+\mu=n$ (this is the K\"{o}nig-Egervary theorem, \cite{LovaPlum}).   A
\textit{K\"{o}nig-Egervary graph} (or simply \textit{KE graph}) is a graph that satisfies this identity. There are non-bipartite KE graphs: the right graph in Figure 2 is an example. KE graphs were  first characterized by Deming  \cite{Deme} and Sterboul \cite{Ster} in 1979. A graph has a \textit{perfect matching} if there is a matching where every vertex of the graph is incident to some edge in the matching (and thus $\mu=\frac{n}{2}$). Deming showed that the problem of determining whether a graph $G$ was KE or not could be reduced to the problem of determining whether a certain extension $G'$ of $G$ with a perfect matching is a KE graph. With respect to a matching $M$, a \textit{blossom} is an odd cycle  where half of one less than the number of edges in the cycle belong to $M$. In this case there must be a unique pair of edges in the cycle which do not belong to $M$. The vertex incident to these two edges is the \textit{blossom tip}. A \textit{blossom pair} is a pair of blossoms whose tips are joined by a path with an odd number of edges, beginning and ending with with edges in $M$ and alternating between edges that are in $M$ and those that are not. Deming proved that if $G$ is a graph with a perfect matching $M$ then, $G$ is a KE graph if, and only if, $G$ contains no blossom pairs. Sterboul gave an equivalent characterization.

Ermelinda DeLaVina's program Graffiti.pc conjectured that, for any graph,  $\alpha=\alpha'$ if, and only if, $\tau=\mu$. This conjecture is proved here for the first time and yields a new characterization of KE graphs.  The Graffiti.pc conjecture can be rewritten: for any graph,  $\alpha=\alpha'$ if, and only if, $\alpha+\mu=n$; or, for a graph $G$, $\alpha(G)=\alpha'(G)$ if, and only if, $G$ is a KE graph. Since 
 a graph was defined to be \emph{totally independence reducible} if $\alpha'=\alpha$, Graffiti.pc's conjecture can also be restated as: a graph $G$ is totally independence reducible if, and only if, $G$ is a KE graph. On the face of it, there is no connection between the non-existence of blossom pairs in a graph and the graph having the property that its independence and critical independence numbers are equal; it is not obvious that this  characterization of KE graphs and the Deming-Sterboul characterization are equivalent.

\begin{thm}
(Graffiti.pc \#329)
For any graph,  $\alpha=\alpha'$ if, and only if, $\tau=\mu$.
\end{thm}

\begin{proof}

Suppose that $\alpha(G)=\alpha'(G)$. It will be shown that $\tau(G)=\mu(G)$ or, equivalently, that $n-\alpha(G)=\mu(G)$.

Let $I$ be a maximum critical independent set. So $\alpha(G)=\alpha'(G)=|I|$. Since $n-\alpha(G)=|N(I)|$, it  remains to show that $\mu(G)=|N(I)|$. Since $I$ is independent, $\mu(G)\leq |N(I)|$. It only remains to show that $\mu(G)\geq |N(I)|$.  Since $I$ is a critical independent set, the Matching Lemma \ref{t:matching} implies that there is a matching from $N(I)$ into $I$ and, thus, that $\mu(G)\geq |N(I)|$.

Suppose now that $\tau(G)=\mu(G)$ or, equivalently, that $n-\alpha(G)=\mu(G)$. It will be shown that $\alpha(G)=\alpha'(G)$.  $\alpha'(G)\leq \alpha(G)$. Suppose $\alpha'(G)<\alpha(G)$. Let $I_c$ be a maximum critical independent set.  Theorem \ref{t:butenko} guarantees the existence of a maximum independent set $J$  such that $I_c\subseteq J$. Since $\mu(G)=n(G)-\alpha(G)$,  $J$ is independent, and $|V\setminus J|=n(G)-\alpha(G)$, there is a matching from $V\setminus J$ into $J$. This implies that each vertex in $N(J)\setminus N(I_c)$ is matched to a vertex in $J\setminus I_c$. So $|J\setminus I_c| \geq |N(J)\setminus N(I_c)|$ .

It will now be shown that $|J| - |N(J)| \geq |I_c| - |N(I_c)|$, implying that $I_c$ is not a maximum critical independent set, as it was assumed to be. $|J| - |N(J)|=(|J\setminus I_c|+|I_c|)-(|N(J)\setminus N(I_c)| + |N(I_c)|) = (|I_c| - |N(I_c))+(|J\setminus I_c| - |N(J)\setminus N(I_c)|) \geq |I_c| - |N(I_c)|$. It follows that $I_c=J$, $|I_c|=|J|$, and $\alpha'(G)=\alpha(G)$, which was to be shown.
\end{proof}

\newpage
 Now Theorem \ref{t:decomposition} can be restated in an interesting and potentially  fruitful way.
 
 \begin{cor}
 For any graph $G$, there is a unique set $X\subseteq V(G)$ such that
\begin{enumerate}
\item $\alpha(G)=\alpha(G[X])+\alpha(G[X^c])$,
\item $G[X]$ is a K\"{o}nig-Egervary graph,
\item for every non-empty independent set $I$ in $G[X^c]$, $|N(I)| > |I|$, and
\item for every maximum critical independent set $J_c$ of $G$, $X=J_c\cup N(J_c)$.
\end{enumerate}
 \end{cor}


\begin{thebibliography}{99}

\bibitem{Agee}
A. Ageev,
On finding  critical independent and vertex sets,
\emph{SIAM J. Discrete Mathematics}, 7:293--295, 1994.




\bibitem{ButeTruk}
S. Butenko and S. Trukhanov,
Using Critical Sets to Solve the Maximum Independent Set Problem,
\textit{Operations Research Letters} 35(4) (2007) 519--524.



\bibitem{Deme}
R. W. Deming,
Independence Numbers of Graphs---an Extension of the Koenig-Egervary Theorem,
\emph{Discrete Mathematics},
27 (1979) 23--33.

\bibitem{GareJohn}
M. Garey and D. Johnson,
Computers and Intractability,
W. H. Freeman and Company, New York, 1979.



\bibitem{Lars}
C. E. Larson,
A Note on Critical Independent Sets,
\emph{Bulletin of the Institute of Combinatorics and its Applications}, v51, Sept. 2007, 34-46.

\bibitem{Lars08}
C. E. Larson,
Graph Theoretic Independence and Critical Independent Sets, dissertation, University of Houston, 2008.

\bibitem{LovaPlum}
L. Lovasz, M. D. Plummer,
Matching Theory,
North Holland, Amsterdam, 1986.

\bibitem{Ster}
F. Sterboul, A characterization of the graphs in which the transversal number equals the matching number, \textit{Journal of Combinatorial Theory. Series B}, 1979 vol:27, 228-229.

\bibitem{Zhan}
C.-Q. Zhang,
Finding critical independent sets and critical vertex subsets are polynomial problems,
\emph{SIAM J. Discrete Mathematics},
3:431--438, 1990.




\end{thebibliography}
\end{document}